\documentclass[11pt,twoside]{amsart}
\usepackage{amsmath}
\usepackage{amsthm}
\usepackage{amsfonts}
\usepackage{amssymb}
\usepackage{latexsym}
\usepackage[all]{xy}

\renewcommand{\uppercasenonmath}[1]{}

\newtheorem{theorem}{Theorem}

\newtheorem*{ack*}{ACKNOWLEDGEMENTS}

\newcommand{\lb}{{\linebreak}}
\newcommand{\dis}{\displaystyle}
\newcommand{\pf}{\noindent\begin {proof}}
\newcommand{\epf}{\end{proof}}

\newcommand{\gre}{\epsilon}

\newcommand{\Hom}{\rm Hom}

\newcommand{\Coker}{\rm Coker }

\newcommand{\bbZ}{\mathbb{bZ}}
\newcommand{\RFlat}{{\rm R \mbox{-}Flat}}
\newcommand{\RProj}{{\rm R \mbox{-}Proj}}
\newcommand{\Rproj}{{\rm R \mbox{-}proj}}
\newcommand{\RMod}{{\rm R \mbox{-}Mod}}

\newcommand{\mcF}{{\mathcal{F}}}

\begin{document}

\footskip30pt

\date{}

\title{Neeman's characterization of $K(\RProj)$ via Bousfield localization}

\author{X.H.\ Fu}
\address{School of Mathematics and Statistics, Northeast Normal University, Changchun, CHINA}
\email{fuxianhui@gmail.com}
\author{Ivo Herzog}

\address{The Ohio State University at Lima, Lima, OH 45804 USA}
\email{herzog.23@osu.edu}

\thanks{The first author is supported by the National Natural Science Foundation of China, Grant No.\ 11301062;
the second by NSF Grant 12-01523.}



\begin{abstract}
Let $R$ be an associative ring with unit and denote by $K(\RProj)$ the homotopy category of complexes of projective left $R$-modules. Neeman proved the theorem that $K(\RProj)$ is $\aleph_1$-compactly generated, with the category $K^+ (\Rproj)$ of left bounded complexes of finitely generated projective $R$-modules providing an essentially small class of such generators. Another proof of Neeman's theorem is explained, using recent ideas of Christensen and Holm, and Emmanouil. The strategy of the proof is to show that every complex in $K(\RProj)$ vanishes in the Bousfield localization $K(\RFlat)/\langle K^+ (\Rproj) \rangle.$       
\end{abstract}

\maketitle

\pagestyle{plain}

Let $R$ be an associative ring with identity and denote by $\RMod$ the category of left modules. The additive subcategories
$$\Rproj \subseteq \RProj \subseteq \RFlat \subseteq \RMod$$
of finitely generated projective, projective and flat modules, respectively, induce full and faithful embeddings 
$$K(\Rproj) \subseteq K(\RProj) \subseteq K(\RFlat) \subseteq K(\RMod)$$
of the associated homotopy categories of complexes, each of which carries the structure of a triangulated category. A.\ Neeman~\cite[Theorem 5.9]{N1} distinguished the homotopy category $K^+ (\Rproj)$ of left bounded complexes of finitely generated projective modules as noteworthy, and if  we let $\langle K^+(\Rproj) \rangle \subseteq K(\RProj)$ be the smallest triangulated subcategory containing $K^+ (\Rproj)$ and closed under direct sums and summands, he established ({\em loc.\ cit.}) the following important equality,
\begin{eqnarray} \label{Eq:N}
K(\RProj) & = &  \langle K^+(\Rproj) \rangle.
\end{eqnarray}
The importance of this equality stems from a preliminary result of Neeman.

\begin{theorem} \label{Neeman} {\cite[Theorem 5.9]{N1}} 
The homotopy category $K^+ (\Rproj)$ of left bounded complexes of finitely generated projective modules is an essentially small subcategory of $\aleph_1$-compact objects of $K(\RFlat).$ 
\end{theorem} 

The main consequence of Equality~(\ref{Eq:N}) is that the category $K(\RProj)$ is an $\aleph_1$-compactly generated triangulated category. This kind of theorem was first proved by P.\ J\o rgensen \cite[Theorem 2.4]{Jo} in the special case of a (left and right) coherent ring over which every flat left $R$-module has finite projective dimension.

In this note, we use a different strategy to prove Equation~(\ref{Eq:N}). Theorem~\ref{Neeman} implies that the the triangulated category $\langle K^+ (\Rproj) \rangle$ is $\aleph_1$-compactly generated, and so satisfies the Brown Representability Theorem (see \cite[Theorem 8.3.3.]{N0}). The inclusion $\langle K^+ (\Rproj) \rangle \subseteq K(\RFlat)$ of triangulated categories therefore~\cite[Proposition 9.1.19]{N0} admits a Bousfield localization. Recall~\cite[Theorem 9.1.13.]{N0} that when such a localization exists, then every complex $P^{\cdot} \in K(\RFlat)$ is part of a distinguished triangle
$${\xymatrix@1@C=45pt
{J^{\cdot} \ar[r]^{u^{\cdot}} & P^{\cdot} \ar[r] & C(u^{\cdot}) \ar[r] & J^{\cdot}[1],
}}$$
where $J^{\cdot} \in \langle K^+(\Rproj) \rangle,$ and the mapping cone $C(u^{\cdot})$ of $u^{\cdot}$ belongs to the right orthogonal subcategory $\langle K^+(\Rproj) \rangle^{\perp}$ in $K(\RFlat),$ defined\footnote{The superscript $\perp$ is to be interpreted throughout relative to the ambient category $K(\RFlat).$} as
$$\langle K^+(\Rproj) \rangle^{\perp} 
:= \{ F^{\cdot} \in K(\RFlat) \; | \; {\Hom}_{K(\RFlat)}(T^{\cdot}, F^{\cdot}) = 0 \;\; \mbox{for all} \;\; 
T^{\cdot} \in \langle K^+(\Rproj) \rangle  \}.$$ 
Emmanouil~\cite{E} provides a comprehensive treatment of the subcategory $\mcF \subseteq K(\RFlat)$ of pure acyclic complexes
(see definition below) and proves a lemma~\cite[Lemma 3.4]{E} that implies 
\begin{eqnarray} \label{Eq:E}
K^+(\Rproj)^{\perp} & \subseteq & \mcF. 
\end{eqnarray}
This clearly implies that $\langle K^+(\Rproj) \rangle^{\perp} \subseteq \mcF.$ \bigskip

To complete the proof of Equality~(\ref{Eq:N}), assuming the inclusion~(\ref{Eq:E}), it suffices to show that if $P^{\cdot} \in K(\RProj),$ then the mapping cone $C(u^{\cdot}) = 0;$ for then $P^{\cdot} \cong J^{\cdot} \in \langle K^+(\Rproj) \rangle.$ As \lb $J^{\cdot} \in \langle K^+(\Rproj) \rangle \subseteq K(\RProj),$ the mapping cone $C(u^{\cdot})$ also belongs to $K(\RProj)$. Let us employ an interesting observation of Christensen and Holm~\cite{CH} to show that 
\begin{eqnarray*} \label{Eq:CH}
K(\RProj) \cap \mcF & = & 0
\end{eqnarray*}
and therefore that $C(u^{\cdot}) = 0.$ 

The argument depends on the result of Benson and Goodearl~\cite{BG} that every periodic flat module is projective. Recall that a module $M$ is {\em periodic} if it is isomorphic to its own syzygy: there is a short exact sequence
$$\vcenter{
\xymatrix@1@C=40pt{
0 \ar[r] & F \ar[r] & P \ar[r] & F \ar[r] & 0,
}}$$
where $P$ is a projective $R$-module. Suppose that $P^{\cdot}: \cdots \to P^{n-1} \to P^n \to P^{n+1} \to \cdots$ is a pure exact complex of projective $R$-modules. Every syzygy $Z^n \subseteq P^n$ is a pure submodule of a projective module, and therefore flat. Consider the complex $\oplus_{n \in \bbZ} \; P^{\cdot}[n]$ associated to $P^{\cdot},$ defined as the coproduct of all the iterated suspensions and desuspensions of $P^{\cdot}.$ It is a repeating complex with the same projective module 
$P = \oplus_{k \in \bbZ} \; P^k$ in every degree and the same {\em flat} syzygy $F = \oplus_{k \in \bbZ} \; Z^k$ in every degree. The flat module $F$ is therefore periodic, and hence projective. Each of its summands $Z^k$ is then projective, so that the complex $P^{\cdot}$ is contractible and $P^{\cdot} = 0$ in $K(\RFlat).$ \bigskip
 
Let us conclude with some remarks on the inclusion~(\ref{Eq:E}). Recall that an exact complex $Y^{\cdot}$ of $R$-modules is \emph{pure acyclic} if the constituent short exact sequences are pure exact. Equivalently, if $C$ is a finitely presented module (concentrated in degree $0$) and integer $n\in \bbZ,$ then ${\Hom}_K(C[-n],Y^{\cdot})=0.$ If  $f : C^{\cdot} \to F^{\cdot}$ is a morphism of complexes where $F^{\cdot}$ is a complex of flat modules and $C^{\cdot}$ a left bounded complex of finitely presented modules, then Emmanouil inductively builds~\cite[Lemma 3.4]{E} a left bounded complex $P^{\cdot}$ of finitely generated modules projective modules such that $f$ factors through $P^{\cdot}.$ This implies, in particular, that if $F^{\cdot}$ belongs to $K^+ (\Rproj)^{\perp} \subseteq K(\RFlat),$ then it must be pure acyclic. His proof does not spell out the initial step of the induction, but this is not a problem; one just backs up by two degrees with a finite projective presentation
$$\xymatrix@C=45pt@R=35pt
{\cdots \ar[r] & 0 \ar[r] \ar@{=}[d] & 0 \ar[r] \ar[d] & 0 \ar[r] \ar[d] & A \ar[r] \ar@{=}[d] & 0 \ar[r] \ar@{=}[d] & \cdots \\
\cdots \ar[r] & 0 \ar[r] \ar[d] & P^{n-2} \ar[r]^{\partial^{n-2}_P} \ar[d]^{f^{n-2}} & P^{n-1} \ar[r]^{\gre^{n-1}} \ar[d]^{f^{n-1}} & A^n \ar[r] \ar[d]^{g^n} & 0 \ar[d] \ar[r] & \cdots \\
\cdots \ar[r] & F^{n-3} \ar[r]^{\partial^{n-3}_F} & F^{n-2} \ar[r]^{\partial^{n-2}_F} & F^{n-1} \ar[r]^{\partial^{n-1}_F} & F^n \ar[r]^{\partial^n_F} & F^{n+1} \ar[r] & \cdots.
}$$
and continues the induction as explained in the proof of~\cite[Lemma 3.2]{E}.

The anonymous referee suggests an illuminating proof of the converse inclusion \lb 
$\mcF \subseteq K(\Rproj)^{\perp} \subseteq K^+ (\Rproj)^{\perp}.$ The idea is to express a complex 
$$\xymatrix@1@C=45pt
{P^{\cdot}=\cdots \ar[r] & P^{n-2} \ar[r]^{\partial_P^{n-2}} & P^{n-1} \ar[r]^{\partial_P^{n-1}} & P^n \ar[r]^{\partial_P^{n}} & P^{n+1} \ar[r] & \cdots
}$$
of finitely generated projective modules as the direct limit of its soft truncations
$$\xymatrix@1@C=45pt
{\sigma_{\leq n}(P^{\cdot})=\cdots \ar[r] & P^{n-2} \ar[r]^{\partial_P^{n-2}} & P^{n-1} \ar[r]^{\partial_P^{n-1}} & P^n \ar[r]^{\gre^n} & A \ar[r] &0\ar[r]& \cdots
}$$ 
with $A = \Coker (\partial_P^{n-1}).$ It is clear that ${\dis P^{\cdot} = \lim_{n \to \infty} \; \sigma_{\leq n}(P^{\cdot}).}$ Moreover, the short exact sequence of chain morphisms
$$\xymatrix{
0 \ar[r] & \bigoplus_{n\in\bbZ}\;  \sigma_{\leq n} (P^{\cdot}) \ar[r] & 
\bigoplus_{n\in\bbZ} \; \sigma_{\leq n} (P^{\cdot}) \ar[r] & {\dis \lim_{n \to \infty} \; \sigma_{\leq n}(P^{\cdot})} = P^{\cdot} \ar[r] & 0
}$$
is chainwise split, which, by~\cite[Proposition 4.15]{I}, yields the triangle
 $$\xymatrix{
P[-1] \ar[r] & \bigoplus_{n\in\bbZ} \; \sigma_{\leq n} (P^{\cdot}) \ar[r] & \bigoplus_{n\in\bbZ} \sigma_{\leq n}(P^{\cdot}) \ar[r] & P}$$ 
in $K(\RMod).$ To prove that a pure acyclic complex $F^{\cdot}$ belongs to $K(\Rproj)^{\perp},$ it suffices to verify that
for every integer $n,$ $\Hom_K (\sigma_{\leq n} (P^{\cdot}), F^{\cdot}) = 0.$ For then
$\Hom_K(\bigoplus_{n\in\bbZ} \; \sigma_{\leq n}(P^{\cdot}), F^{\cdot}) =
\prod_{n\in\bbZ}\; \Hom_K(\sigma_{\leq n} (P^{\cdot}), F^{\cdot}) = 0.$

Consider the truncation
$$\xymatrix@1@C=45pt{
\tau_{\leq n}(P^{\cdot}) = \cdots \ar[r] & P^{n-2} \ar[r]^{\partial_P^{n-2}} & P^{n-1} \ar[r]^{\partial_P^{n-1}} & P^n \ar[r]^{\partial_P^{n}} & 0 \ar[r] & \cdots.
}$$
The obvious short exact sequence of chain morphisms 
$$\xymatrix@1@C=45pt{
0 \ar[r] & A[-n-1] \ar[r] & \sigma_{\leq n} (P^{\cdot}) \ar[r] & \tau_{\leq n} (P^{\cdot}) \ar[r] & 0
}$$
induces a distinguished triangle in $K(\RMod).$ Since $F^{\cdot}$ is pure acyclic and each $\tau_{\leq n} (P^{\cdot})$ is a right bounded complex of projective modules, it follows that $\Hom_K (A[-n], F^{\cdot}) = \Hom_K (\tau_{\leq n} (P^{\cdot}), F^{\cdot}) = 0.$ By applying the functor $\Hom_K (-, F^{\cdot})$ to the induced triangle one obtains that 
$\Hom_K (\sigma_{\leq n} (P^{\cdot}), F^{\cdot}) = 0,$ as required. 
\bigskip

Together, these elegant ideas of Emmanouil and the referee provide us with a new proof of Neeman's equality~\cite[Theorem 8.6]{N1} 
$$K(\RProj)^{\perp}= \mcF.$$
Clearly, $K(\RProj)^{\perp} \subseteq K^+({\Rproj})^{\perp} = \mcF.$ On the other hand, $K^+(\Rproj)^{\perp} \subseteq \mcF$ implies that $K(\RProj)^{\perp} = \langle K^+ (\Rproj) \rangle^{\perp} \subseteq \mcF,$ by Equality~(\ref{Eq:N}). \bigskip

\noindent {\em Acknowledgements.} The suggestion above of the referee resulted in a substantial improvement of our original argument. In fact, the initial submission of this note contained several proofs that have since appeared in publication; we thank S.\ Estrada and S.\ Iyengar for pointing these out. We are also grateful to J.\ \v{S}\'{t}ov\'{i}\v{c}ek for important comments on an earlier draft.

\end{document}